\documentclass[11pt, letterpaper]{article}
\usepackage[utf8]{inputenc}
\usepackage{amssymb, amsfonts, amsmath, amsthm, mathtools}
\usepackage[margin=1.2in]{geometry}
\usepackage{enumitem}
\usepackage{ifthen}
\usepackage[none]{hyphenat}
\usepackage{theoremref}
\usepackage[usenames,dvipsnames]{xcolor}
\usepackage{tikz}
\usepackage{color}

\newtheoremstyle{mystyle}
{1.2em} 
{0.1em} 
    {\itshape} 
    {} 
    {\bfseries} 
    {.} 
    {.5em} 
    {} 
\theoremstyle{mystyle}
\newtheorem{thm}{Theorem}[section]
\theoremstyle{mystyle}
\newtheorem{prop}[thm]{Proposition}
\theoremstyle{mystyle}

\theoremstyle{mystyle}
\newtheorem{obs}[thm]{Observation}
\theoremstyle{mystyle}

\theoremstyle{mystyle}
\newtheorem{cor}[thm]{Corollary}
\theoremstyle{mystyle}

\theoremstyle{mystyle}
\newtheorem{que}{Question}
\providecommand{\boksie}{\ensuremath{\mathbin{\raisebox{0.3mm}{$\scriptstyle\square$}}}}
\input{tcilatex}

\definecolor{lgray}{gray}{0.95}
\definecolor{mgray}{gray}{0.40}
\usetikzlibrary{fit,positioning,calc}
\tikzstyle{std}=[ circle, draw=black,fill=black, inner sep=0pt, minimum size=2mm]
\tikzstyle{bred}=[circle, draw=black,fill=red,thick,  inner sep=2pt, minimum size=2mm]
\tikzstyle{bgreen}=[ circle, draw=black,fill=green,thick,  inner sep=2pt, minimum size=2.5mm]
\tikzstyle{sqRed}=[rectangle, draw=black,fill=red,thick,  inner sep=2pt, minimum size=2.5mm]
\tikzstyle{trEdge}=[color=black, line width = 1pt, style=dotted]

\begin{document}

\title{\textbf{Total Roman Domination Edge-Critical Graphs}}
\author{C. Lampman\thanks{%
Supported by an Undergraduate Student Research Award from the Natural
Sciences and Engineering Research Council of Canada.}, C. M. Mynhardt\thanks{%
Supported by a Discovery Grant from the Natural Sciences and Engineering
Research Council of Canada.}, S. E. A. Ogden\thanks{%
Supported by a Science Undergraduate Research Award from the University of
Victoria.} \\
Department of Mathematics and Statistics\\
University of Victoria\\
Victoria, BC, \textsc{Canada}\\
{\small chloe.lampman3@gmail.com, kieka@uvic.ca, sogden@uvic.ca}}
\date{}
\maketitle

\begin{abstract}
A total Roman dominating function on a graph $G$ is a function $%
f:V(G)\rightarrow \{0,1,2\}$ such that every vertex $v$ with $f(v)=0$ is
adjacent to some vertex $u$ with $f(u)=2$, and the subgraph of $G$ induced
by the set of all vertices $w$ such that $f(w)>0$ has no isolated vertices.
The weight of $f$ is $\Sigma _{v\in V(G)}f(v)$. The total Roman domination
number $\gamma _{tR}(G)$ is the minimum weight of a total Roman dominating
function on $G$. A graph $G$ is $k$-$\gamma _{tR}$-edge-critical if $\gamma
_{tR}(G+e)<\gamma _{tR}(G)=k$ for every edge $e\in E(\overline{G})\neq
\emptyset $, and $k$-$\gamma _{tR}$-edge-supercritical if it is $k$-$\gamma
_{tR}$-edge-critical and $\gamma _{tR}(G+e)=\gamma _{tR}(G)-2$ for every
edge $e\in E(\overline{G})\neq \emptyset $. We present some basic results on $\gamma_{tR}$-edge-critical graphs and characterize certain classes of $\gamma
_{tR}$-edge-critical graphs. In addition, we show that, when $k$ is small, there is a connection between $k$-$\gamma _{tR}$-edge-critical graphs and graphs which are critical with respect to the domination and total domination numbers.      
\end{abstract}



\section{Introduction}


We consider the behaviour of the total Roman domination number of a graph $G$ upon the addition of edges to $G$. A \emph{dominating set} $S$ in a graph $G$ is a set of vertices such that
every vertex in $V(G)-S$ is adjacent to at least one vertex in $S$.
The \emph{domination number} $\gamma (G)$ is the cardinality of a minimum
dominating set in $G$. A \emph{total dominating set} $S$ (abbreviated by 
\emph{TD-set}) in a graph $G$ with no isolated vertices is a set of vertices
such that every vertex in $V(G)$ is adjacent to at least one vertex in $S$.
The \emph{total domination number} $\gamma _{t}(G)$ (abbreviated by \emph{%
TD-number}) is the cardinality of a minimum total dominating set in $G$. For $S\subseteq V(G)$ and a function $f:S\rightarrow \mathbb{R}$, define $%
f(S)=\Sigma _{s\in S}f(s)$. A \emph{Roman dominating function} (abbreviated
by \emph{RD-function}) on a graph $G$ is a function $f:V(G)\rightarrow
\{0,1,2\}$ such that every vertex $v$ with $f(v)=0$ is adjacent to some
vertex $u$ with $f(u)=2$. The \emph{weight} of $f$, $\omega (f)$, is defined
as $f(V(G))$. The \emph{Roman domination number} $\gamma _{R}(G)$
(abbreviated by \emph{RD-number}) is defined as $\min \{\omega (f):f\text{
is an RD-function on }G\}$. For an RD-function $f$, let $V_{f}^{i}=\{v\in
V(G):f(v)=i\}$ and $V_{f}^{+}=V_{f}^{1}\cup V_{f}^{2}$. Thus, we can uniquely express an RD-function $f$ as $f=(V_f^0,V_f^1,V_f^2)$. 

As defined by Ahanger, Henning, Samodivkin and Yero \cite{AHSY}, a \emph{%
total Roman dominating function} (abbreviated by \emph{TRD-function}) on a
graph $G$ with no isolated vertices is a Roman dominating function with the
additional condition that $G[V_{f}^{+}]$ has no isolated vertices. The \emph{%
total Roman domination number} $\gamma _{tR}(G)$ (abbreviated by \emph{%
TRD-number}) is the minimum weight of a TRD-function on $G$, that is, $%
\gamma _{tR}(G)=\min \{\omega (f):f\text{ is a TRD-function on }G\}$. A
TRD-function $f$ such that $\omega (f)=\gamma _{tR}(G)$ is called a $\gamma
_{tR}(G)$-\emph{function}, or a $\gamma _{tR}$-\emph{function} if the graph $%
G$ is clear from the context; $\gamma_{R}$\emph{-functions} are defined analogously.

The addition of an edge to a graph has the potential to change its total
domination or Roman domination number. Van der Merwe, Mynhardt and Haynes 
\cite{MMH} studied $\gamma _{t}$-\emph{edge-critical graphs}, that is,
graphs $G$ for which $\gamma _{t}(G+e)<\gamma _{t}(G)$ for each $e\in E(%
\overline{G})$ and $E(\overline{G})\neq\emptyset $. We consider the same
concept for total Roman domination. A graph $G$ is \emph{total Roman
domination edge-critical}, or simply $\gamma_{tR}$\emph{-edge-critical}, if $%
\gamma _{tR}(G+e)<\gamma _{tR}(G)$ for every edge $e\in E(\overline{G})$ and 
$E(\overline{G})\neq \emptyset $. We say that $G$ is $k$\emph{-}$\gamma _{tR}$\emph{%
-edge-critical} if $\gamma _{tR}(G)=k$ and $G$ is $\gamma _{tR}$%
-edge-critical. If $\gamma _{tR}(G+e)\leq \gamma_{tR}(G)-2$ for every edge $%
e\in E(\overline{G})$ and $E(\overline{G})\neq\emptyset $, we say that $G$
is $\gamma _{tR}$\emph{-edge-supercritical}. If $\gamma_{tR}(G+e)=\gamma _{tR}(G)$ for all $e\in E(\overline{G})$, or $E(\overline{G%
})=\emptyset $, we say that $G$ is \emph{stable}.

Pushpam and Padmapriea \cite{PP} established bounds on the total Roman domination number of a graph in terms of its order and girth. Total Roman domination in trees was studied by Amjadi, Nazari-Moghaddam, Sheikholeslami and Volkmann \cite{ANSV}, as well as by Amjadi,  Sheikholeslami and Soroudi \cite{ASS1}. The authors of \cite{ASS1} also studied Nordhaus-Gaddum bounds for total Roman domination in \cite{ASS2}. Campanelli and Kuziak \cite{CK} considered total Roman domination in the lexicographic product of graphs. We refer the reader to the well-known books \cite{CL} and \cite{HHS} for
graph theory concepts not defined here. Frequently used or lesser known
concepts are defined where needed.

We begin with some general results regarding the addition of an edge $e\in E(%
\overline{G})$ to a graph $G$ in Section \ref{Sec:Adding an edge}. In
Section \ref{Sec:ntR}, we characterize $n$-$\gamma _{tR}$-edge-critical
graphs of order $n$. We characterize $4$-$\gamma _{tR}$%
-edge-critical graphs in Section \ref{Sec:4tR}, and, after investigating $\gamma_{tR}$-edge-supercritical graphs in Section \ref{Sec:Super}, we present a necessary condition for $5$-$\gamma _{tR}$-edge-critical graphs in Section \ref{Sec:5tR}. In Section \ref{Sec:Spider}, we
determine the total Roman domination number of spiders and characterize $%
\gamma _{tR}$-edge-critical spiders. As can be expected, every graph $G$ with $\gamma_{tR}(G)=k\geq4$ is
a spanning subgraph of a $k$-$\gamma _{tR}(G)$-edge-critical graph; a short
proof is given in Section \ref{Sec:Diameter}, where we also show that for
any $k\geq 4$, there exists a $k$-$\gamma _{tR}$-edge-critical graph of
diameter $2$. We conclude in Section \ref{Sec:Future} with ideas for future
research.

\section{Adding an edge}

\label{Sec:Adding an edge}

We begin with a result from \cite{HMM} which bounds the effect the addition
of an edge can have on the total domination number of a graph and show that
the same bounds hold with respect to the total Roman domination number.

\begin{prop}
\thlabel{Myn1} \emph{\cite{HMM}} For a graph $G$ with no isolated vertices, if $uv\in E(\overline{G})$, then $\gamma_{t}(G)-2\leq \gamma _{t}(G+uv)\leq \gamma _{t}(G)$.
\end{prop}

\vspace{1.2mm}

An edge $uv\in E(\overline{G})$ is \emph{critical} if $\gamma
_{tR}(G+uv)<\gamma _{tR}(G)$. The following proposition restricts the
possible values assigned to the vertices of a critical edge $uv$ by a $%
\gamma _{tR}(G+uv)$-function $f$, which will be useful in proving subsequent results. For a graph $G$ and a vertex $v\in V(G)$, the \emph{open neighbourhood} of $v$ in $G$ is $N_{G}(v)=\{u\in
V(G):uv\in E(G)\}$, and the \emph{closed neighbourhood} of $v$ in $G$ is $%
N_{G}[v]=N_{G}(v)\cup \{v\}$. When $G\neq K_{2}$, the unique neighbour of an
end-vertex of $G$ is called a \emph{support vertex}.  

\begin{prop}
\thlabel{set added edge}
Given a graph $G$ with no isolated vertices, if $uv\in E(\overline{G})$ is a critical edge and $f$
is a $\gamma _{tR}(G+uv)$-function, then $\{f(u),f(v)\}\in
\{\{2,2\},\{2,1\},\{2,0\},\{1,1\}\}$. If, in addition, $\deg (u)=\deg (v)=1$%
, then there exists a $\gamma_{tR}(G+uv)$-function $f$ such that $f(u)=f(v)=1$.
\end{prop} 

\noindent \emph{Proof.} Let $G$ be a graph with no isolated vertices, $uv\in E(\overline{G})$ such
that $\gamma _{tR}(G+uv)<\gamma _{tR}(G)$, and $f$ a $\gamma _{tR}$-function
on $G+uv$. Suppose for a contradiction that $\{f(u),f(v)\}\notin
\{\{2,2\},\{2,1\},\{2,0\},\{1,1\}\}$. Then $\{f(u),f(v)\}\in
\{\{0,0\},\{0,1\}\}$. Note that, in either case, the edge $uv$ cannot affect
whether $u$ and $v$ are dominated, or whether, in the case where (say) $%
f(v)=1$, $v$ is isolated. Hence $f$ is a TRD-function of $G$, contradicting $\gamma _{tR}(G+uv)<\gamma _{tR}(G)$. Therefore $%
\{f(u),f(v)\}\in \{\{2,2\},\{2,1\},\{2,0\},\{1,1\}\}$.

Now, suppose in addition that $\deg (u)=\deg (v)=1$, and let $f$ be a $%
\gamma _{tR}(G+uv)$-function such that $|V_{f}^{2}|$ is as small as
possible. Let $w$ and $x$ be the unique neighbours of $u$ and $v$,
respectively, noting that possibly $w=x$. Suppose for a contradiction that $%
f(u)=2$ (without loss of generality). If $f(v)=0$, then $f(w)>0$, otherwise $u$ would be isolated in $G[V_f^+]$. Thus, regardless of whether $w=x$ or not, consider the function $f^{\prime
}:V(G)\rightarrow \{0,1,2\}$ defined by $f^{\prime }(u)=f^{\prime }(v)=1$
and $f^{\prime }(y)=f(y)$ for all other $y\in V(G)$. Otherwise, if $%
f(v)\geq 1$, then clearly $f(w)=0$. Thus, regardless of whether $w=x$ or not, consider the function $f^{\prime }:V(G)\rightarrow \{0,1,2\}$ defined by $f^{\prime }(u)=f^{\prime
}(w)=1$ and $f^{\prime }(y)=f(y)$ for all other $y\in V(G)$. In either case, 
$f^{\prime }$ is a $\gamma _{tR}$-function on $G+uv$. However, $%
|V_{f^{\prime }}^{2}|<|V_{f}^{2}|$, contradicting $|V_{f}^{2}|$ being as
small as possible. Hence $f(u)\neq 2$, and thus $f(u)=f(v)=1$.~$\square $

\begin{prop}
\thlabel{tR bounds} 
Given a graph $G$ with no isolated vertices, if ${uv\in E(\overline{G})}$, then $\gamma
_{tR}(G)-2\leq \gamma _{tR}(G+uv)\leq \gamma _{tR}(G)$.
\end{prop}

\noindent\emph{Proof.} Let $G$ be a graph with no isolated vertices. Clearly, adding an edge cannot increase the total Roman domination number, hence the upper bound holds. Now, let $uv\in E(\overline{G})$. Note that when $\gamma_{tR}(G+uv) =\gamma_{tR}(G)$ the lower bound clearly holds. So assume $\gamma_{tR}(G+uv)
< \gamma_{tR}(G)$ and let $f$ be a $\gamma_{tR}(G+uv)$-function. By 
\thref{set added edge}, $\{f(u), f(v)\} \in \{ \{2,2\}, \{2,1\}, \{2,0\}, \{1,1\} \}$.

First assume $\{f(u),f(v)\}\in \{\{2,2\},\{2,1\},\{1,1\}\}$. Then $f$ is a
RD-function of $G$, and the only possible isolated vertices in $G[V_{f}^{+}]$
are $u$ and $v$. Consider the function $f^{\prime }:V(G)\rightarrow
\{0,1,2\} $ defined as follows: If $u$ is isolated in $G[V_{f}^{+}]$, choose 
$u^{\prime }\in N_{G}(u)$ and let $f^{\prime }(u^{\prime })=1$. Similarly,
if $v$ is isolated in $G[V_{f}^{+}]$, choose $v^{\prime }\in N_{G}(v)$ and
let $f^{\prime }(v^{\prime })=1$. Let $f^{\prime }(x)=f(x)$ for all other $x\in
V(G)$. Now, assume instead that $f(u)=2$ and $f(v)=0$ (without loss of generality). Since $%
u$ is not isolated in $G[V_{f}^{+}]$, $f$ is a TRD-function of $G-v$.
Consider the function $f^{\prime }:V(G)\rightarrow \{0,1,2\}$ defined as
follows: Let $f^{\prime }(v)=1$. Then, if $v$ is isolated in $%
G[V_{f^{\prime }}^{+}]$, choose $v^{\prime }\in N_{G}(v)$ and let $f^{\prime
}(v^{\prime })=1$. Let $f^{\prime }(x)=f(x)$ for all other $x\in V(G)$. In either case, $f^{\prime }$ is a TRD-function of $G$ and $\omega(f^{\prime })\leq \gamma _{tR}(G+uv)+2$. Thus $\gamma _{tR}(G)\leq \gamma
_{tR}(G+uv)+2$, and hence the lower bound holds.~$\square $

\section{$\protect\gamma _{tR}$-Edge-critical graphs with large TRD-numbers}

\label{Sec:ntR}

We now investigate the $\gamma _{tR}$-edge-critical graphs $G$ which have the
largest TRD-number, namely $|V(G)|$. A \emph{subdivided star} is a tree obtained from a star on at least three
vertices by subdividing each edge exactly once. A \emph{double star} is a
tree obtained from two disjoint non-trivial stars by joining the two central
vertices (choosing either central vertex in the case of $K_{2}$). The \emph{%
corona} $\func{cor}(G)$ (sometimes denoted by $G\circ K_{1}$) of $G$ is
obtained by joining each vertex of $G$ to a new end-vertex.

Connected graphs $G$ for which $\gamma _{tR}(G)=|V(G)|$ were characterized
in \cite{AHSY}. There, Ahanger et al. defined $\mathcal{G}$ as the family of
connected graphs obtained from a $4$-cycle $v_{1},v_{2},v_{3},v_{4},v_{1}$
by adding $k_{1}+k_{2}\geq 1$ vertex-disjoint paths $P_{2}$, and joining $%
v_{i}$ to the end of $k_{i}$ such paths, for $i\in \{1,2\}$. Note that
possibly $k_{1}=0$ or $k_{2}=0$. Furthermore, they defined $\mathcal{H}$ to
be the family of graphs obtained from a double star by subdividing each
pendant edge once and the non-pendant edge $r\geq 0$ times. For $r\geq 0$%
, define $\mathcal{H}_{r}\subseteq \mathcal{H}$ as the family of graphs in $%
\mathcal{H}$ where the non-pendant edge was subdivided $r$ times. 

\begin{prop}
\thlabel{Hen1} \emph{\cite{AHSY}} If $G$ is a connected graph of order $%
n\geq 2$, then $\gamma _{tR}(G)=n$ if and only if one of the following
holds. 
\vspace{-2mm}
\setlist{nolistsep}
\begin{enumerate}
\item[$(i)$] $G$ is a path or a cycle.

\item[$(ii)$] $G$ is the corona of a graph.

\item[$(iii)$] $G$ is a subdivided star.

\item[$(iv)$] $G\in \mathcal{G}\cup \mathcal{H}$.
\end{enumerate}
\end{prop}

\vspace{1.2mm}

Using \thref{Hen1}, we characterize connected $n$-$\gamma _{tR}$-edge-critical graphs
as follows.

\begin{thm}
\thlabel{n edge-crit} A connected graph $G$ of order $n\geq 4$ is $n$-$\gamma _{tR}$%
-edge-critical if and only if $G$ is one of the following graphs: %
\vspace{-2mm}
\setlist{nolistsep}
\begin{enumerate}
\item[$(i)$] $C_{n}$, $n\geq 4$,

\item[$(ii)$] $\func{cor}(K_{r})$, $r\geq 3$,

\item[$(iii)$] a subdivided star of order $n\geq 7$,

\item[$(iv)$] $G\in \mathcal{G}$,

\item[$(v)$] $G\in \mathcal{H}-\mathcal{H}_{0}-\mathcal{H}_{2}$.
\end{enumerate}
\end{thm}

\noindent \emph{Proof.} Let $G$ be a connected graph of order $n\geq4$ with $\gamma
_{tR}(G)=n$. First, suppose $G$ is any of the graphs listed in $(i)-(v)$
above. Then, for any $e\in E(\overline{G})$, $G+e$ is not one of the graphs
listed in \thref{Hen1}. Therefore $\gamma _{tR}(G+e)<n$ for all $e\in E(%
\overline{G})$, and thus $G$ is $\gamma _{tR}$-edge-critical.

Otherwise, suppose $G$ is not one of the graphs listed in $(i)-(v)$ above.
Note that since $\gamma _{tR}(G)=n$, $G$ is still listed in \thref{Hen1} $%
(i)-(iv)$. If $G\cong P_{n}:v_{1},...,v_{n}$, $n\geq4$, then $%
G+v_{1}v_{n}\cong C_{n}$ and $\gamma _{tR}(G)=\gamma _{tR}(C_{n})=n$. If $%
G\cong \func{cor}(F)$, where $F$ is not a complete graph of order at least $3$, then $\gamma_{tR}(G)=%
\gamma _{tR}(G+uv)$ for any $uv\in E(\overline{F})$. If $G$ is a subdivided
star of order less than $7$, then $G=P_{5}$. In each of these cases, $G$ is
clearly not $\gamma _{tR}$-edge-critical.

Now consider $G\in \mathcal{H}$. Let $w_{1},...,w_{k}$ be the leaves of $G$, 
$u_{1},...,u_{k}$ be their respective support vertices, and $v_{1},...,v_{m}$
be the path such that $v_{1}$ and $v_{m}$ are the two support vertices in
the original double star $S$, labelled so that $w_{1}$ is adjacent, in $S$,
to $v_{1}$. Note that $m=r+2$, and therefore $m\geq 2$. If $G\in \mathcal{H}%
_{0}$, consider the graph $G+v_{2}w_{1}$, and note that $G+v_{2}w_{1}\in 
\mathcal{G}$. Therefore, by \thref{Hen1}, $\gamma _{tR}(G+v_{2}w_{1})=n$,
and thus $G$ is not $\gamma _{tR}$-edge-critical. Similarly, if $G\in 
\mathcal{H}_{2}$, consider the graph $G+v_{1}v_{4}$, and note that $%
G+v_{1}v_{4}\in \mathcal{G}$. Therefore, by \thref{Hen1}, $%
\gamma_{tR}(G+v_{1}v_{4})=n$, and again $G$ is not $\gamma _{tR}$%
-edge-critical.~$\square $

\section{$4$-$\protect\gamma _{tR}$-Edge-critical graphs}

\label{Sec:4tR}

Before we characterize the graphs $G$ such that $\gamma _{tR}(G)=4$ and $%
\gamma _{tR}(G+e)=3$ for any $e\in E(\overline{G})$ (that is, the graphs
which are $4$-$\gamma _{tR}$-edge-critical), we present the following result
from \cite{PP} which characterizes the graphs with a total Roman domination
number of $3$, the smallest possible TRD-number. Note that while the authors
required that $G$ has girth $3$, the result actually holds in general for
any graph $G$ on at least $3$ vertices, as we now show. A \emph{universal
vertex }of $G$ is a vertex that is adjacent to all other vertices of $G$.

\begin{prop}
\thlabel{tR=3} For a graph $G$ of order $n\geq 3$ with no isolated vertices, $\gamma _{tR}(G)=3$ if and
only if $\Delta (G)=n-1$, that is, $G$ has a universal vertex.
\end{prop}

\noindent \emph{Proof.} Suppose $\gamma _{tR}(G)=3$ and let $%
f=(V_{f}^{0},V_{f}^{1},V_{f}^{2})$ be a $\gamma _{tR}(G)$-function.
If $V_{f}^{2}=\emptyset $, then $|V_{f}^{1}|=3$, and thus $n=3$.
Since $G$ has no isolated vertices, this implies that $G=K_{3}$ or $P_{3}$,
both of which have a universal vertex. Otherwise, assume $|V_{f}^{2}|=1$ and 
$|V_{f}^{1}|=1$. Pick $u,v\in V(G)$ so that $f(u)=1$ and $f(v)=2$. Since $%
G[V_{f}^{+}]$ has no isolated vertices, $uv\in E(G)$. Furthermore, since $%
\gamma _{tR}(G)=3$, $f(x)=0$ for all other $x\in V(G)$. Therefore 
$N_G[v]=V(G)$, and thus $v$ is a universal vertex.

Conversely, suppose $G$ has a universal vertex $v$, and take any $u\in
N_{G}(v)$. Consider the TRD-function $f:V(G)\rightarrow \{0,1,2\}$ defined by $%
f(v)=2$, $f(u)=1$, and $f(x)=0$ for all other $x\in V(G)$. Since $G$ has at least three vertices, $\gamma_{tR}(G)>2$. Therefore, since $\omega (f)=3$, we conclude that $\gamma_{tR}(G)=3$.~$\square $

\vspace{1.2mm}

A \emph{galaxy} is defined as the disjoint union of two or more non-trivial
stars. The characterization of $4$-$\gamma _{tR}$-edge-critical graphs
follows; note that this class of graphs is exactly the class of $2$-$\gamma$%
-edge-critical graphs, as characterized by Sumner and Blitch \cite{SB}.

\begin{thm}
\thlabel{4-tR edge-crit} A graph $G$ with no isolated vertices is $4$-$\gamma _{tR}$-edge-critical if
and only if $\overline{G}$ is a galaxy.
\end{thm}

\noindent \emph{Proof.} Let $G$ be a graph of order $n$ with no isolated vertices. Suppose first that $G$ is $4$-$\gamma _{tR}$-edge-critical. Then for any $e\in E(\overline{G})$, $\gamma _{tR}(G+e)=3$, and thus \thref{tR=3} implies that the addition of any edge to $G$ creates
a universal vertex. Therefore, for each edge $uv\in E(\overline{G})$, one
of $u$ and $v$ has degree $n-2$ in $G$; that is, one of $u$ and $v$ is a leaf in $%
\overline{G}$. Since each edge of $\overline{G}$ connects a leaf to either a support vertex or another leaf, the components of $\overline{G}$ are non-trivial stars.
Moreover, $\overline{G}$ has at least two components, otherwise $G$ has an
isolated vertex. 

Conversely, suppose $\overline{G}$ is a galaxy. Since $%
\overline{G}$ has no isolated vertices, $G$ has no universal vertices, and thus, by %
\thref{tR=3}, $\gamma _{tR}(G)>3$. Let $u$ and $v$ be vertices in different
components of $\overline{G}$, and define $f:V(G)\rightarrow \{0,1,2\}$ by $f(u)=f(v)=2$ and $f(x)=0$ for all other $x\in V(G)$. Clearly $f$ is a TRD-function on $G$, and hence $\gamma _{tR}(G)=4$. Since the
deletion of any edge in $\overline{G}$ produces an isolated vertex, the
addition of any edge to $G$ creates a universal vertex. Therefore, by %
\thref{tR=3}, $\gamma _{tR}(G+e)=3$ for all $e\in E(\overline{G})$, and
hence $G$ is $4$-$\gamma _{tR}$-edge-critical.~$\square $

\begin{cor}
If $G$ is a connected $(n-2)$-regular graph, then $G$ is $4$-$\gamma _{tR}$%
-edge-critical.
\end{cor}

\vspace{1.2mm}

Having characterized $4$-$\gamma_{tR}$-edge-critical graphs, our
next result demonstrates the existence of stable graphs with total Roman
domination number $4$.

\begin{prop}
\thlabel{n-3 reg} If $G$ is an $(n-3)$-regular graph of order $n\geq 6$,
then $\gamma _{tR}(G)=4$. Moreover, $G$ is stable.
\end{prop}

\noindent \emph{Proof.} We prove that $\gamma (G)=2$. Since $G$ is $(n-3)$%
-regular, its complement $\overline{G}$ is $2$-regular. If $\overline{G}$ is
disconnected, let $u$ and $v$ be vertices in different components of $%
\overline{G}$. Otherwise, if $\overline{G}$ is connected, then $\overline{G}\cong C_{n}$%
, $n\geq 6$, and thus we can choose $u,v\in V(\overline{G})$ such that $d_{\overline{G}%
}(u,v)\geq 3$. In either case, $N_{\overline{G}}[u]\cap N_{\overline{G}%
}[v]=\emptyset $. In $G$, $u$ dominates all vertices in $G-N_{\overline{G}%
}(u)$ and $v$ dominates all vertices in $G-N_{\overline{G}}(v)$. Therefore $\{u,v\}$ dominates $G$, and thus, since $G$ has no universal vertex, $\gamma (G)=2$.

Now, define $f:V(G)\rightarrow \{0,1,2\}$ by $f(u)=f(v)=2$ and $f(y)=0$ for
all other $y\in V(G)$. Since $uv\in E(G)$, $f$ is a TRD-function
on $G$ and $\omega (f)=4$, so $\gamma _{tR}(G)\leq 4$. Since $G$ has no
universal vertex, $\gamma _{tR}(G)>3$ by \thref{tR=3}, and thus $\gamma
_{tR}(G)=4$, as required. Furthermore, since the addition of any edge to $G$
does not create a universal vertex, it follows from \thref{tR=3} that $\gamma
_{tR}(G+e)=\gamma _{tR}(G)$ for all $e\in E(\overline{G})$. Therefore $G$ is
stable.~$\square $

\section{$\protect\gamma _{tR}$-Edge-supercritical graphs}

\label{Sec:Super}

We now consider the graphs $G$ which attain the lower bound in 
\thref{tR bounds} for all $e\in E(\overline{G})$, that is, $\gamma _{tR}$%
-edge-supercritical graphs. An edge $uv\in E(\overline{G})$ is \emph{supercritical} if $\gamma_{tR}(G+uv)=\gamma_{tR}(G)-2$. Haynes, Mynhardt and Van der Merwe \cite{HMM}
defined a graph $G$ to be $\gamma _{t}$\emph{-edge-supercritical} if $%
\gamma_{t}(G+e)=\gamma _{t}(G)-2$ for all $e\in E(\overline{G})$. We begin
with their characterization of $\gamma _{t}$-edge-supercritical graphs.

\begin{prop}
\thlabel{Myn2} \emph{\cite{HMM}} A graph $G$ is $\gamma _{t}$%
-edge-supercritical if and only if $G$ is the union of two or more
non-trivial complete graphs.
\end{prop}

\vspace{1.2mm}

The proof of the previous result relies on the fact that, if $u$ and $v$ are vertices of a graph $G$ with $d(u,v)=2$, then $\gamma_t(G)-1\leq \gamma_t(G+uv)$. However, the analogous result does not hold with respect to the total Roman domination number, as we now show. Consider the graph $G=\func{cor}(K_3)$. By \thref{Hen1}, $\gamma_{tR}(G)=6$. Consider any two non-adjacent vertices $u$ and $v$ in $G$ such that $\deg(u)=1$ and $\deg(v)=3$. Clearly $uv$ is a supercritical edge with $d(u,v)=2$, and thus $d(u,v)=2$ does not always imply that $\gamma_{tR}(G)-1\leq \gamma_{tR}(G+uv)$.

As a result, the classification of $\gamma_{tR}$-edge-supercritical graphs will be less straightforward than that of $\gamma_{t}$-edge-supercritical graphs. However, it is easy to see that there are no $5$-$\gamma_{tR}$-edge-supercritical graphs, the smallest possible TRD-number of a $\gamma_{tR}$-edge-supercritical graph, and that the disjoint union of two or more complete graphs of order at least $3$ is $\gamma_{tR}$-edge-supercritical.  

\begin{prop} 
\thlabel{no 5-super}
\setlist{nolistsep}
\begin{enumerate}
\item[]
\item[$(i)$] There are no $5$-$\gamma_{tR}$-edge-supercritical graphs.
\item[$(ii)$] If $G$ is the disjoint union of $k\geq2$ complete graphs, each of order at least $3$, then $G$ is $3k$-$\gamma_{tR}$-edge-supercritical. 
\end{enumerate}
\end{prop}  
\vspace{-1mm}

\noindent \emph{Proof.}
\vspace{-3mm} 
\setlist{nolistsep}
\begin{enumerate}
\item[$(i)$] Suppose for a contradiction that $G$ is a $5$-$\gamma_{tR}$-edge-supercritical graph. Then $\gamma(G+uv)=3$ for any edge $uv\in E(\overline{G})$. However, as in the proof of \thref{4-tR edge-crit}, this implies that $\overline{G}$ is a galaxy, that is, $G$ is $4$-$\gamma_{tR}$-edge-critical, a contradiction.
\item[$(ii)$] It follows from \thref{tR=3} that $\gamma_{tR}(G)=3k$. Moreover, joining any two vertices in different components of $G$ results in a graph with TRD-number $3k-2$.~$\square $
\end{enumerate}

\section{$5$-$\protect\gamma _{tR}$-Edge-critical graphs}

\label{Sec:5tR}

We now investigate the graphs which are $5$-$%
\gamma _{tR}$-edge-critical. We begin with the following results from \cite%
{AHSY}, which bound $\gamma _{tR}(G)$ in terms of $\gamma _{t}(G)$.

\begin{prop}
\thlabel{Hen2} \emph{\cite{AHSY}} If $G$ is a graph with no isolated vertices, then $\gamma _{t}(G)\leq \gamma _{tR}(G)\leq 2\gamma _{t}(G)$. Furthermore, $\gamma _{tR}(G)=\gamma _{t}(G)$ if and only if $G$ is the disjoint union of copies of $K_2$. 
\end{prop}

\vspace{1.2mm}

Note that Amjadi et al. \cite{ANSV} characterized the trees which attain the upper bound in \thref{Hen2}.

\begin{prop}
\thlabel{Hen3} \emph{\cite{AHSY}} Let $G$ be a connected graph of order $n\geq3$. Then $\gamma _{tR}(G)=\gamma _{t}(G)+1$ if and only if $\Delta(G)=n-1$, that is, $G$ has a universal vertex.    
\end{prop}

\vspace{1.2mm}

By \thref{tR=3}, \thref{Hen3} implies that, if $G$ is a connected graph
of order $n\geq 3$, then $\gamma _{tR}(G)=\gamma _{t}(G)+1$ if and only if $\gamma _{tR}(G)=3$. These results lead to the following observation.

\begin{obs}
\thlabel{tR>t+1} If $G$ is a connected graph of order $n\geq 3$ such that $%
\Delta (G)\leq n-2$, then $\gamma _{t}(G)+2\leq \gamma _{tR}(G)\leq 2\gamma
_{t}(G)$.
\end{obs}

\vspace{1.2mm}

We now provide a result characterizing graphs with $\gamma _{tR}\in \{3,4\}$
in terms of their domination and total domination numbers that will be
useful in describing $5$-$\gamma _{tR}$-edge-critical graphs.

\begin{prop}
\thlabel{t=2 iff tR=34} If $G$ is a connected graph of order $n\geq 3$, then 
$\gamma _{tR}(G)\in \{3,4\}$ if and only if $\gamma _{t}(G)=2$. Moreover, $\gamma (G)=1$ when $\gamma _{tR}(G)=3$, and $\gamma (G)=2$ when $%
\gamma _{tR}(G)=4$.
\end{prop}

\noindent \emph{Proof.} Suppose first that $\gamma _{t}(G)=2$. By \thref{Hen2}, $2\leq \gamma _{tR}(G)\leq 4$. Clearly $\gamma _{tR}(G)\neq 2$, since $n\geq 3$. Therefore $\gamma _{tR}(G)\in \{3,4\}$.

Conversely, suppose $\gamma _{tR}(G)\in\{3,4\}$. First, if $\gamma _{tR}(G)=3
$, then \thref{tR=3} implies that $G$ has a universal vertex. Therefore $%
\gamma_{t}(G)=2$ and $\gamma(G)=1$. Otherwise, if $\gamma _{tR}(G)=4$, then
\thref{tR=3} implies that $G$ has no universal vertex. Therefore, by \thref{tR>t+1}, $%
\gamma_{t}(G)+2 \leq 4$, and thus $\gamma_{t}(G)=2$. Furthermore, since $\gamma(G)\leq\gamma_t(G)$ and $G$
has no universal vertex, $\gamma(G)=2$.~$\square $

\begin{prop}
\thlabel{5-tR edge-crit} For any graph $G$, if $G$ is $5$-$\gamma _{tR}$-edge-critical, then $G$ is either $3$-$\gamma _{t}$-edge-critical or $G=K_{2}\cup
K_{n}$ for $n\geq 3$, in which case $G$ is $4$-$\gamma _{t}$%
-edge-supercritical.
\end{prop}

\noindent \emph{Proof.} Suppose $G$ is $5$-$\gamma _{tR}$-edge-critical. By %
\thref{t=2 iff tR=34}, $\gamma _{t}(G)>2$ and $\gamma _{t}(G+e)=2$ for any $%
e\in E(\overline{G})$. Therefore, by \thref{Myn1}, $G$ is either $3$-$\gamma
_{t}$-edge-critical or $4$-$\gamma _{t}$-edge-supercritical. If $G$ is $4$-$%
\gamma _{t}$-edge-supercritical, then by \thref{Myn2}, $G$ is the union of
two or more non-trivial complete graphs. Since $\gamma _{tR}(G)=5$, this implies that $%
G=K_{2}\cup K_{n}$ for $n\geq 3$.~$\square $

\vspace{1.2mm}

Note that if we add the condition that $G$ is not $6$-$\gamma_{tR}$-edge-supercritical, then the above becomes a necessary and sufficient condition. Clearly $G=K_{2}\cup K_{n}$ is $5$-$%
\gamma _{tR}$-edge-critical for any $n\geq 3$. Otherwise, if $G$ is $3$-$\gamma _{t}$-edge-critical, then by \thref{t=2 iff tR=34}, $\gamma_{tR}(G)>4$ and $\gamma_{tR}(G+e)\in\{3,4\}$ for any $e\in E(\overline{G})$. By \thref{Hen2}, $\gamma_{tR}(G)\leq 6$, and thus, since $G$ is not $6$-$\gamma_{tR}$-edge-supercritical, $\gamma_{tR}(G)=5$. Hence $G$ is $5$-$\gamma_{tR}$-edge-critical, as required.

\section{$\protect\gamma _{tR}$-Edge-critical spiders}

\label{Sec:Spider}A \emph{(generalized) spider }$\mathrm{Sp}%
(l_{1},...,l_{k}),\ l_{i}\geq 1,\,k\geq 2$, is a tree obtained from the star 
$K_{1,k}$ with centre $u$ and leaves $v_{1},...,v_{k}$ by subdividing the
edge $uv_{i}$ exactly $l_{i}-1$ times, $i=1,...,k$. Thus, a spider $\mathrm{Sp}%
(2,...,2)$ is a subdivided star. The $u-v_{i}$ paths (of length $l_{i}$) are
called the \emph{legs} of the spider, while $u$ is its \emph{head}. We now
investigate the spiders which are $\gamma _{tR}$-edge-critical. Note that when $k=2$, $\mathrm{Sp}(l_{1},...,l_{k})\cong P_n$ for $n\geq3$, which, by \thref{n edge-crit}, is not $\gamma _{tR}$-edge-critical. We begin with two propositions restricting $\gamma _{tR}$-edge-criticality in general
graphs, which will be useful in classifying $\gamma _{tR}$-edge-critical
spiders.

\begin{prop}
\thlabel{end deg 3} For a graph $G$ with no isolated vertices, if $G$ has an end-vertex $w$ with support vertex $x$ such that $G[N(x)-\{w\}]$ is not complete, then $G$ is not $\gamma _{tR}$-edge-critical.
\end{prop}

\noindent \emph{Proof.} Suppose $u,v\in N_{G}(x)-\{w\}$ such that $uv\in E(%
\overline{G})$. We claim that $\gamma _{tR}(G+uv)=\gamma _{tR}(G)$. Suppose
for a contradiction that $\gamma _{tR}(G+uv)<\gamma _{tR}(G)$, and consider
a $\gamma _{tR}$-function $f=(V_{f}^{0},V_{f}^{1},V_{f}^{2})$ on $G+uv$.
Note that, since $w$ is an end-vertex, $f(x)>0$. By \thref{set added edge}, $%
\{f(u),f(v)\}\in \{\{2,2\},\{2,1\},\{2,0\},\{1,1\}\}$. Since $ux,vx\in E(G)$
and at least one of $f(u)$ and $f(v)$ is positive, we can assume without loss of
generality that $f(x)=2$. In any case, $f$ is also a TRD-function on $G$,
contradicting $\gamma _{tR}(G+uv)<\gamma _{tR}(G)$. Therefore $\gamma
_{tR}(G+uv)=\gamma _{tR}(G)$ and $G$ is not $\gamma _{tR}$-edge-critical.~$%
\square $

In a tree, the support vertex of a leaf is called a \emph{stem}. A stem is
called \emph{weak} if it is adjacent to exactly one leaf, and \emph{strong}
if it is adjacent to two or more leaves. A vertex $b$ of a tree such that $%
\deg (b)\geq 3$ is called a \emph{branch vertex}. An \emph{endpath }in a
tree is a path from a branch vertex to a leaf, all of whose internal
vertices have degree $2$.\emph{\ } The next result follows immediately from %
\thref{end deg 3}.

\begin{cor}
\thlabel{strong stems} If $T$ is a $\gamma _{tR}$-edge-critical tree, then $%
T $ contains no stems of degree at least $3$, and hence no strong stems.
\end{cor}

\begin{prop}
\thlabel{no long legs} For a graph $G$ with no isolated vertices, if $G$ has two endpaths $%
v_{0}, v_{1}, ..., v_{k}$ and $u_{0}, u_{1}, ..., u_{m}$, where $k,m\geq 3$ and $v_{k}$
and $u_{m}$ are leaves, then $G$ is not $\gamma _{tR}$-edge-critical.
\end{prop}

\noindent \emph{Proof.} We claim that $\gamma _{tR}(G+v_{k}u_{m})=\gamma
_{tR}(G)$. Suppose for a contradiction that $\gamma
_{tR}(G+v_{k}u_{m})<\gamma _{tR}(G)$, and let $f$ be a $\gamma _{tR}$%
-function on $G+v_{k}u_{m}$. Then, by \thref{set added edge}, we may assume $%
f(u_{m})=f(v_{k})=1$. Define $f^{\prime }:V(G)\rightarrow \{0,1,2\}$ as
follows: If $f(v_{k-1})=0$, then clearly $f(v_{k-2})=2$ and $f(v_{k-3})\geq
1 $, so let $f^{\prime }(v_{k-1})=f^{\prime }(v_{k-2})=1$. Otherwise, let $f^{\prime }(v_{k-1})=f(v_{k-1})$ and $f^{\prime }(v_{k-2})=f(v_{k-2})$. Similarly, if $%
f(u_{m-1})=0$, then clearly $f(u_{m-2})=2$ and $f(u_{m-3})\geq 1$, so let $%
f^{\prime }(u_{m-1})=f^{\prime }(u_{m-2})=1$. Otherwise, let $f^{\prime }(u_{m-1})=f(u_{m-1})$ and $f^{\prime }(u_{m-2})=f(u_{m-2})$. Finally, let $f^{\prime
}(y)=f(y)$ for all other $y\in V(G)$. Clearly $f^{\prime }$ is a
TRD-function on $G$ and $\omega (f^{\prime })=\omega (f)$, contradicting $%
\gamma _{tR}(G+v_{k}u_{m})<\gamma _{tR}(G)$. Therefore $\gamma
_{tR}(G+v_{k}u_{m})=\gamma _{tR}(G)$, and thus $G$ is not $\gamma _{tR}$%
-edge-critical.~$\square $

\vspace{1.2mm}

Let $S$ be a spider with $k\geq 3$ legs. In what follows, let $c$ be the
head of $S$, and let the $k$ legs be labelled $c,v_{i1},v_{i2},...,v_{im_{i}}
$, where $i\in \{1,2,...,k\}$, in order of increasing length. Let $m=m_{k}$,
that is, $m$ is the length of a longest leg of $S$. We begin by determining
the TRD-number of spiders.

\begin{prop}
\thlabel{tR spider} If $S$ is a spider of order $n$ with $k\geq 3$ legs such
that $S$ has $y$ legs of length $2$, then 
\begin{equation*}
\gamma _{tR}(S)=%
\begin{cases}
n & \text{if $y\geq k-1$} \\ 
n-k+y+1 & \text{if $1\leq y<k-1$} \\ 
n-k+2 & \text{if $y=0.$}%
\end{cases}%
\end{equation*}
\end{prop}

\vspace{-1mm}
\noindent \emph{Proof.} Suppose $S$ has $x$ legs of length $1$, and consider
a $\gamma _{tR}$-function $f$ on $S$ such that $|V_{f}^{2}|$ is as small as
possible. First, suppose $y\geq k-1$. If $y=k$, then $S$ is a subdivided
star. Otherwise, if $y=k-1$, then $S$ has exactly one leg that is not of
length $2$, and thus either $x=1$ or $x=0$. If $x=1$, then $S$ is the corona
of a graph (specifically, $S=\func{cor}(K_{1,k-1})$). Otherwise, if $x=0$,
then $m=m_{k}\geq 3$, and $S\in \mathcal{H}_{r}$, where $r=m-3$. In any
case, by \thref{Hen1}, $\gamma _{tR}(S)=n$.

Assume therefore that $y<k-1$. Then $S$ has at least two legs that are not
of length $2$. Therefore $S$ is not one of the graphs listed in \thref{Hen1}%
, and thus $\gamma _{tR}(S)<n$. Hence there is some vertex $u\in V(S)$ such
that $f(u)=2$ and $f(w)=0$ for at least two vertices $w$ adjacent to $u$.
Furthermore, since $f$ is a TRD-function, such a vertex $u$ is not isolated
in $S[V_{f}^{+}]$, and thus $\deg (u)\geq 3$. Since $c$ is the only vertex in 
$S$ with degree at least $3$, $f(c)=2$. Therefore $c$ Roman dominates each $%
v_{i1}$, and we need $f(v_{i1})$ to be positive for at least one $i$ to
ensure that $S[V_{f}^{+}]$ has no isolated vertices.

Consider an arbitrary leg $c,v_{i1},v_{i2},...,v_{im_{i}}$ of $S$. If $%
m_{i}=1$, then $f(v_{i1})\in \{0,1\}$ in order for $f$ to totally Roman
dominate $c$ and $v_{i1}$. If $m_{i}=2$, a total weight of $2$ on $v_{i1}$
and $v_{i2}$ is required in order for $f$ to total Roman dominate $%
\{v_{i1},v_{i2}\}$. Since $|V_{f}^{2}|$ is as small as possible, $%
f(v_{i1})=f(v_{i2})=1$. Finally, if $m_{i}>2$, by \thref{Hen1} and since $%
f(c)=2$, a total weight of at least $m_{i}-1$ on $v_{i1},...,v_{im_{i}}$ is
required in order for $f$ to totally Roman dominate $c$ and $%
\{v_{i1},...,v_{im_{i}}\}$. Moreover, by the choice of $f$, $f(v_{i1})\in
\{0,1\}$ and $f(v_{i2})=\cdots =f(v_{im})=1$. Therefore $\omega (f)\geq
n-k+y+1$.

Now, if $y>0$, where (say) $m_{j}=2$, then $f(v_{j1})=1$. By minimality and
since $c$ is adjacent to $v_{j1}$, $f(v_{i1})=0$ for each $i$ such that $%
m_{i}\neq 2$. Then $\gamma _{tR}(S)=\omega (f)=n-k+y+1$, as required.
Otherwise, if $y=0$, then $f(v_{i1})=1$ for some $i$ to ensure that $c$ is
not isolated in $S[V_{f}^{+}]$. By minimality, $f(v_{j1})=0$ for each $j\neq
i$. Therefore $\gamma _{tR}(S)=\omega (f)=n-k+2$.~$\square $

\vspace{1.2mm}

The characterization of $\gamma _{tR}$-edge-critical spiders follows. Our
result also shows that a spider of order $n$ is $\gamma _{tR}$-edge-critical
if and only if it is $n$-$\gamma _{tR}$-edge-critical.

\begin{thm}
\thlabel{edge-crit spider} A spider $S=\func{Sp}(l_{1},...,l_{k}),$\ $k\geq 3
$, is $\gamma _{tR}$-edge-critical if and only if $l_{i}=2$ for each $i$, $%
1\leq i\leq k-1$, and $l_{k}\in \{2,m\}$, where $m=4$ or $m\geq 6$.
\end{thm}

\noindent \emph{Proof.} Suppose $S$ has order $n$. If $l_{i}=2$ for each $%
i=1,...,k$, then $S$ is a subdivided star and, by \thref{n edge-crit}, $S$
is $n$-$\gamma _{tR}$-edge-critical. Now, suppose $S$ has exactly one leg of
length $m\neq 2$. If $m=1$, then by \thref{end deg 3}, $S$ is not $\gamma
_{tR}$-edge-critical. If $m=3$ or $m=5$, then $S\in \mathcal{H}_{r}$ with $%
r=0$ or $2$, respectively, and thus, by \thref{n edge-crit}, $S$ is not $%
\gamma _{tR}$-edge-critical. If $m=4$ or $m\geq 6$, then $S\in \mathcal{H}%
_{r}$ with $r=m-3$, and therefore, by \thref{n edge-crit}, $S$ is $n$-$%
\gamma _{tR}$-edge-critical. Finally, suppose $S$ has at least two legs that
are not of length $2$. Again, by \thref{end deg 3}, if $S$ has a leg of
length $1$, $S$ is not $\gamma _{tR}$-edge-critical. Otherwise, assume $S$
has at least two legs of length at least $3$. Then, by \thref{no long legs}, 
$S$ is not $\gamma _{tR}$-edge-critical.~$\square $

\section{$k$-$\protect\gamma _{tR}$-Edge-critical graphs with minimum
diameter}

\label{Sec:Diameter}

We now consider the minimum diameter possible in a $k$-$%
\gamma _{tR}$-edge-critical graph, for $k\geq 4$. There are no $\gamma _{tR}$%
-edge-critical graphs with diameter $1$, as the only graphs with diameter $1$
are non-trivial complete graphs, which are clearly not $\gamma _{tR}$-edge-critical
since $E(\overline{G})=\emptyset $. Therefore, the minimum possible diameter
for a $\gamma _{tR}$-edge-critical graph is $2$. Asplund, Loizeaux and Van der Merwe \cite{ALM} constructed families of $3$-$\gamma_{t}$-edge-critical graphs with diameter $2$. We will show that, for any $%
k\geq 4$, there exists a $k$-$\gamma _{tR}$-edge-critical graph of diameter $%
2$. We first present the following proposition which shows that every graph $%
G$ without a dominating vertex is a spanning subgraph of a $\gamma _{tR}(G)$%
-edge-critical graph with the same total Roman domination number, which will
be useful in proving our result.

\begin{prop}
\thlabel{span} For a graph $G$ with no isolated vertices, if ${\gamma _{tR}(G)=k\geq 4}$, then $G$ is a spanning subgraph of a $k$-$\gamma _{tR}(G)$-edge-critical graph.
\end{prop}

\noindent \emph{Proof.} Suppose $\gamma _{tR}(G)=k\geq4$. If $G$ is $k$-$%
\gamma _{tR}(G)$-edge-critical, then we are done. Otherwise, there is, by
definition, some edge $e_{1}\in E(\overline{G})$ such that $\gamma
_{tR}(G+e_{1})=\gamma _{tR}(G)$. Let $G_{1}=G+e_{1}$. If $G_{1}$ is $k$-$%
\gamma _{tR}(G)$-edge-critical, then we are done. Otherwise, there is some
edge $e_{2}\in E(\overline{G_{1}})$ such that $\gamma
_{tR}(G_{1}+e_{2})=\gamma _{tR}(G_{1})$. Let $G_{2}=G_{1}+e_{2}$. Continuing
in this way, we eventually obtain a graph $G_{i}$ such that for all $e\in E(%
\overline{G_{i}})$, $\gamma _{tR}(G_{i}+e)<\gamma _{tR}(G_{i})$ and $\gamma
_{tR}(G_{i})=\gamma _{tR}(G_{i-1})=\cdots =\gamma _{tR}(G_{1})=\gamma
_{tR}(G)$. Since $k\geq 4$, $G_{i}$ is not complete and thus $E(G_{i})\neq
\emptyset $. Therefore, $G_{i}$ is a $k$-$\gamma _{tR}(G)$-edge-critical
graph, of which $G$ is a spanning subgraph.~$\square $

\vspace{1.2mm}

Before demonstrating the existence of $k$-$\gamma _{tR}$-edge-critical
graphs of diameter $2$ for any $k\geq4$, we determine the TRD-number of $%
K_{n}\boksie K_{m}$, where $n,m\geq 2$. Consider the vertices of $K_{n}%
\boksie K_{m}$ as an $n\times m$ matrix, where vertices $v_{ij}$ and $v_{st}$
are adjacent if and only if $i=s$ or $j=t$. The rows and columns of the
matrix form disjoint copies of $K_{m}$ and $K_{n}$, respectively. If $v_{ij}$
and $v_{st}$ are nonadjacent, then $v_{sj}$ is adjacent to both $v_{ij}$ and 
$v_{st}$, and hence $\func{diam}(K_{n}\boksie K_{m})=2$. 

\begin{prop}
\thlabel{tR(KnXKm)} If $m$ and $n$ are integers such that $m\geq n\geq 2$,
then $\gamma _{tR}(K_{n}\boksie K_{m})=2n$.
\end{prop}

\noindent \emph{Proof.} Let $G=K_{n}\boksie K_{m}$. To see that $\gamma
_{tR}(G)\leq 2n$, consider the TRD-function $g=(V_{g}^{0},V_{g}^{1},V_{g}^{2})$ on $G$ where $V_{g}^{1}=\emptyset $ and $V_{g}^{2}=\{v_{i1}:1\leq i\leq n\}$.

Now, suppose for a contradiction that $\gamma _{tR}(G)\leq 2n-1$ and
consider a TRD-function $f=(V_{f}^{0},V_{f}^{1},V_{f}^{2})$ on $G$ with $%
\omega (f)=2n-1$. Each vertex $v$ dominates one row and one column of $G$,
so if $|V_{f}^{2}|=x$ (note that $x\leq n-1$), then at most $x$ rows and at
most $x$ columns are dominated by elements of $V_{f}^{2}$. Let $S$ be the
set of vertices undominated by $V_{f}^{2}$. Then $|S|\geq (n-x)(m-x)\geq
(n-x)^{2}$. Moreover, $|V_{f}^{1}|=(2n-1)-2x$ since $\omega (f)=2n-1$ and $|V_{f}^{2}|=x$.

If $x=n-1$, then $|V_{f}^{1}|=1$. Since $f$ is a TRD-function and $%
|S|\geq(n-x)^{2}$, $|S|=1$; say $S=\{w\}$. Hence $V_{f}^{1}=\{w\}$. However, $%
V_{f}^{2}$ does not dominate $w$, and thus $w$ is isolated in $G[V_{f}^{+}]$%
, which is a contradiction. Therefore, there is no TRD-function on $G$ with
weight $2n-1$ when $x=n-1$.

\setlength{\abovedisplayskip}{3pt}
\setlength{\belowdisplayskip}{3pt}
Otherwise, if $x<n-1$, then
\begin{align*}
|S|-|V_{f}^{1}|& \geq (n-x)^{2}-(2n-1-2x) \\
& =x^{2}-2(n-1)x+(n-1)^{2} \\
& =(n-1-x)^{2} \\
& >0.
\end{align*}%
Therefore, the number of vertices undominated by $V_{f}^{2}$ is greater than 
$|V_{f}^{1}|$, contradicting $f$ being a TRD-function. Thus there is no
TRD-function on $G$ with weight $2n-1$ when $x<n-1$. We conclude that $%
\gamma _{tR}(G)=2n$.~$\square $

\begin{thm}
\thlabel{tR diam 2} If $k\geq 4$, then there exists a $k$-$\gamma _{tR}$%
-edge-critical graph of diameter~$2$.
\end{thm}

\noindent \emph{Proof.} First, assume that $k$ is even; say $k=2l$ for some $%
l\geq 2$. Let $G_{l}=K_{l}\boksie K_{l}$. By \thref{tR(KnXKm)}, $\gamma
_{tR}(G_{l})=2l$, and, by \thref{span}, $G_{l}$ is a spanning subgraph of a $%
k $-$\gamma _{tR}$-edge-critical graph $G_{l}^{\prime }$. Since $k>3$, \thref{tR=3} implies that $G_{l}^{\prime }$ has no dominating vertex, and hence $2\leq 
\func{diam}(G_{l}^{\prime })\leq \func{diam}(G_{l})=2$.

Now, consider the case where $k$ is odd; say $k=2l+1$ for some $l\geq 2$.
Let $G_{l}^{d}$ be the graph formed by taking $K_{l+1}\boksie K_{l+1}$ and
deleting the vertices in the set $\{v_{j1}:\lfloor \frac{l}{2}\rfloor +2\leq
j\leq l+1\}$. Similarly to $G_{l}$, $\func{diam}(G_{l}^{d})=2$. See Figure $1$.

We claim that $\gamma_{tR}(G_l^d)=2l+1$. To see that $\gamma_{tR}(G_l^d)\leq
2l+1$, consider the following TRD-function on $G_l^d$: If $l$ is even, place
two $2$'s in each of the first $\frac{l}{2}-1$ rows, and one $2$ in each of
rows $\frac{l}{2}$ and $\frac{l}{2}+1$, such that they span columns $2$
through $l+1$. At this point, every vertex in $G_l^d$ is dominated. However,
the $2$'s in rows $\frac{l}{2}$ and $\frac{l}{2}+1$ are isolated, so place a $1$
in row $\frac{l}{2}$ such that it shares a column with the $2$
in row $\frac{l}{2}+1$. Otherwise, if $l$ is odd, place two $2$'s in each of
the first $\frac{l-1}{2}$ rows, and one $2$ in row $\frac{l+1}{2}$, such
that they span columns $2$ through $l+1$. Similarly to the even case, every
vertex in $G_l^d$ is now dominated. However, the $2$ in row $\frac{l+1}{2}$
is isolated, so place a $1$ in row $\frac{l-1}{2}$ such that it shares a
column with that $2$. In either case, we have a TRD-function on $G_l^d$ with
weight $2l+1$, hence $\gamma_{tR}(G_l^d)\leq 2l+1$.

Now, suppose for a contradiction that $\gamma _{tR}(G_{l}^{d})<2l+1$, and
consider a TRD-function $f=(V_{f}^{0},V_{f}^{1},V_{f}^{2})$ on $G_{l}^{d}$
with $\omega (f)=2l$. We claim that $f(v_{j1})=0$ for all $1\leq j\leq
\lfloor \frac{l}{2}\rfloor +1$. If $f(v_{j1})=2$ for $x\geq 1$ vertices in
column $1$, the undominated vertices in columns $2$ through $l+1$ form the
graph $K_{l}\boksie K_{l+1-x}$. By \thref{tR(KnXKm)}, a TRD-function on $%
K_{l}\boksie K_{l+1-x}$ requires a weight of $2\min \{l,l+1-x\}=2(l+1-x)$.
However, since $2x+2(l+1-x)>2l$, this is impossible. Therefore $%
f(v_{j1})\neq 2$ for all $1\leq j\leq \lfloor \frac{l}{2}\rfloor +1$. If $%
f(v_{j1})=1$ for $x\geq 1$ vertices in column $1$, the undominated vertices
in columns $2$ through $l+1$ (that is, those for which $f$ could be assigned
a $2$) form the graph $K_{l}\boksie K_{l+1}$. Again by \thref{tR(KnXKm)}, a
TRD-function on $K_{l}\boksie K_{l+1}$ requires a weight of $2\min
\{l,l+1\}=2l$. However, $x+2l>2l$ for $x\geq 1$, so this is also not
possible. Therefore, $f(v_{j1})=0$ for all $1\leq j\leq \lfloor \frac{l}{2}%
\rfloor +1$.

{\begin{figure}
	\centering
	\begin{tikzpicture}[scale=1.2]	

		\begin{scope}[shift={(-2.5,0)}] 
			\foreach \j in {1,2,3}
				\foreach \i in {1,2,3}
					{\ifthenelse {\i = 1}{\node [std] (v\i\j) at (\j, -\i)[label=above: $2$]{};}{\node [std] (v\i\j) at (\j, -\i){};}}

			\foreach \i in {1,2,3}
				{\draw (v\i1)--(v\i3);
				\draw (v1\i)--(v3\i);
				\draw (v\i1) to [out=45,in=135] (v\i3);
				\draw (v1\i) to [out=225,in=135] (v3\i);}
		\end{scope}

		\begin{scope}[shift={(2,0.5)}] 
			\foreach \j in {1,2,3,4}
				\foreach \i in {1,2}
					{\ifthenelse {\(\j=2 \OR \j=3\) \AND \i=1}{\node [std] (v\i\j) at (\j, -\i)[label=above: $2$]{};}{\ifthenelse {\j=4 \AND \i=2}{\node [std] (v\i\j) at (\j, -\i)[label=right: $2$]{};}}{\ifthenelse {\j=4 \AND \i=1}{\node [std] (v\i\j) at (\j, -\i)[label=above: $1$]{};}{\node [std] (v\i\j) at (\j, -\i){};}}}
					
			\foreach \j in {2,3,4}
				\foreach \i in {3,4} 
					\node [std] (v\i\j) at (\j, -\i){};

			\draw (v11)--(v21);

			\foreach \i in {1,2}
				{\draw (v\i1)--(v\i4);
				\draw (v\i1) to [out=45,in=135] (v\i3);
				\draw (v\i2) to [out=45,in=135] (v\i4);
				\draw (v\i1) to [out=45,in=135] (v\i4);}

			\foreach \i in {3,4}
				{\draw (v\i2)--(v\i4);
				\draw (v\i2) to [out=45,in=135] (v\i4);}

			\foreach \j in {2,3,4}
				{\draw (v1\j)--(v4\j);
				\draw (v1\j) to [out=225,in=135] (v3\j);
				\draw (v2\j) to [out=225,in=135] (v4\j);
				\draw (v1\j) to [out=225,in=135] (v4\j);}

		\end{scope}
	\end{tikzpicture}	
	\caption{The graphs $G_3$ and $G_3^d$ with a $\gamma_{tR}$-function}
\end{figure}
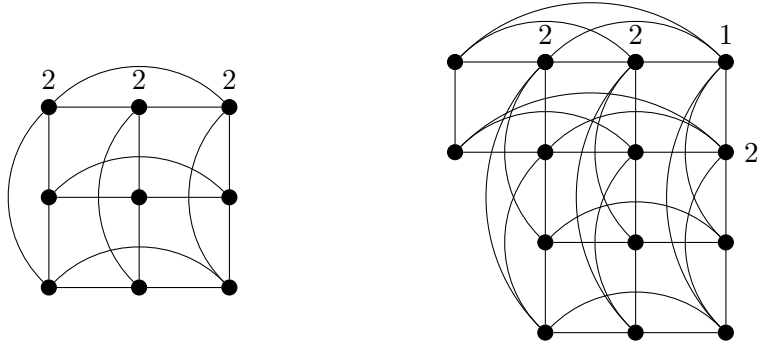}

As a result, in order to totally Roman dominate the first column, there must
be a $2$ in each of the first $\lfloor \frac{l}{2}\rfloor +1$ rows, none of
which can be in the first column. That is, for each $1\leq s\leq \lfloor 
\frac{l}{2}\rfloor +1$, $f(v_{st})=2$ for some $2\leq t\leq l+1$. Let $S$ be
the set of these vertices. Note that, thus far, we have accounted for a
total weight of 
\begin{equation*}
2\left( \left\lfloor \frac{l}{2}\right\rfloor +1\right) =%
\begin{cases}
l+2 & \text{if }l\text{ is even} \\ 
l+1 & \text{if }l\text{ is odd,}%
\end{cases}%
\end{equation*}%
which leaves a weight of $l-2$ if $l$ is even and $l-1$ if $l$ is odd to be assigned. That is, a weight of $2(\lceil \frac{l}{2}\rceil -1)$ remains to be accounted for. We now claim that no two vertices in $S$ can be in the same column. If the vertices in $S$ span fewer than $\lfloor \frac{l}{2}\rfloor +1$ columns, then the vertices which are undominated by $S$ induce a graph containing $K_{\lceil \frac{l}{2}\rceil }\boksie K_{\lceil \frac{l}{2}\rceil }$ as subgraph. If $l=2$, then no weight remains to dominate this vertex, as $2(\lceil \frac{l}{2}\rceil -1)=0$. Otherwise, if $l>2$, \thref{tR(KnXKm)} implies that $\gamma _{tR}(K_{\lceil \frac{l}{2}\rceil }\boksie K_{\lceil \frac{l}{2}\rceil })= 2(\lceil \frac{l}{2}\rceil )$. However, $2(\lceil \frac{l}{2}\rceil )> 2(\lceil \frac{l}{2}\rceil -1)$. In either case, this contradicts $f$ being a TRD-function, and thus no vertices of $S$ share a column.

Therefore, the vertices left undominated by $S$ induce a graph $T\cong
K_{\lceil \frac{l}{2}\rceil }\boksie K_{\lceil \frac{l}{2}\rceil -1}$, with $%
\lceil \frac{l}{2}\rceil $ rows and $\lceil \frac{l}{2}\rceil -1$ columns.
Moreover, the vertices in $S$ are all isolated, as none share a row or
column. By \thref{tR(KnXKm)}, $\gamma _{tR}(T)=2(\lceil \frac{l}{2}\rceil
-1) $. Thus the entire remaining weight is required in order to dominate $T$;
necessarily, the vertices in $V_{f}^{+}-S$ belong to rows and columns that
do not contain vertices in $S$. However, this still leaves the vertices in $%
S $ isolated, which contradicts $f$ being a TRD-function on $G_{l}^{d}$.
Therefore $\gamma _{tR}(G_{l}^{d})\geq 2l+1$ and we conclude that $\gamma
_{tR}(G_{l}^{d})=2l+1$. As in the case where $k$ is even, $G_{l}^{d}$ is a
spanning subgraph of a $k$-$\gamma _{tR}$-edge-critical graph with diameter~$%
2$.~$\square $

\section{Future work}

\label{Sec:Future}

We showed in Section \ref{Sec:Super} that the disjoint union of two of more complete graphs, each having order at least $3$, is $\gamma_{tR}$-edge-supercritical. We also explained that a proof similar to that of \thref{Myn2} does not work for total Roman domination. Hence we pose the following question.

\begin{que}
Are the disjoint unions of two or more complete graphs, each having order at least $3$, the only $\gamma_{tR}$-edge-supercritical graphs?
\end{que}

\vspace{1.2mm}

Note that if this is the case, \thref{5-tR edge-crit} automatically becomes a necessary and sufficient condition for a graph to be $5$-$\gamma_{tR}$-edge-critical.  

Now consider, for a moment, Roman dominating functions, and
suppose a graph $G$ has non-adjacent vertices $u$ and $v$ such that $%
f(u)=f(v)=0$ for every $\gamma _{R}$-function $f$ on $G$. We claim that $%
\gamma _{R}(G+uv)=\gamma _{R}(G)$. Suppose $\gamma _{R}(G+uv)<\gamma _{R}(G)$
and let $f$ be a $\gamma _{R}$-function on $G+uv$. Similar to 
\thref{set added edge}, we may assume without loss of generality that $f(u)=2$ and $%
f(v)=0$, otherwise $f$ is a RD-function on $G$ such that $\omega (f)<\gamma
_{R}(G)$. However, the function $f^{\prime }$ defined by $f^{\prime
}(v)=1$ and $f^{\prime }(y)=f(y)$ for all other $y\in V(G)$ is a $%
\gamma _{R}$-function on $G$ such that $f^{\prime }(v)>0$, contrary to our
assumption. The situation for total Roman domination is different.

For a graph $G$, we define $u\in V(G)$ to be a \emph{dead vertex} if every $%
\gamma _{tR}$-function $f$ on $G$ has $f(u)=0$. Not only do there exist
graphs $G$ containing non-adjacent dead vertices $u$ and $v$ such that $%
\gamma _{tR}(G+uv)<\gamma _{tR}(G)$, but it is possible to find such a graph 
$G$ with $\gamma _{tR}(G+uw)<\gamma _{tR}(G)$ for every edge $uw\in E(%
\overline{G})$, that is, every edge in $E(\overline{G})$ incident with the
dead vertex $u$ is critical. We define the graph $D_{n}$ below and show that 
$D_{n}$ is such a graph.

Let $D_{n}$ be the graph composed of $n\geq 2$ copies of $K_{4}-e$ sharing a
single central vertex as follows: Let $c$ be the central vertex, $%
w_{1},...,w_{n}$ be the degree two vertices, and $u_{1},...,u_{n}$ and $%
v_{1},...,v_{n}$ be the remaining vertices (where $u_{i}$ and $v_{i}$ are
adjacent for each $i$) such that $c,u_{i},w_{i},v_{i},c$ is a $4$-cycle in $%
D_{n}$ for each $1\leq i\leq n$. See Figure $2$. 

\begin{prop}
\thlabel{tR(D_n)} If $n\geq 2$, then $\gamma _{tR}(D_{n})=2n+1$. Moreover, $%
w_{i}$ is a dead vertex for each $1\leq i\leq n$.
\end{prop}

\noindent \emph{Proof.} To see that $\gamma _{tR}(D_{n})\leq 2n+1$, consider the TRD-function $g:V(D_{n})\rightarrow \{0,1,2\}$ on $D_n$ defined by $g(c)=1$, $g(u_{i})=2$ for $%
1\leq i\leq n$, and $g(y)=0$ for all other $y\in V(D_{n})$.

We claim that, if $f$ is a TRD-function on $D_{n}$ with $\omega (f)\leq
2n+1 $, then $f(c)=1$. If $f(c)=2$, then the only vertices that remain
undominated in $D_{n}$ are $w_{i}$ for $1\leq i\leq n$. However, since $%
d(w_{i},w_{j})=4$ for all $i\neq j$, a weight of $2n$ is required in order
to totally Roman dominate these vertices, contradicting $\omega (f)\leq 2n+1$%
. If $f(c)=0$, then since $D_{n}-c$ is the disjoint union of $n$ triangles,
\thref{Hen1} implies that a weight of $3n$ is required to totally Roman dominate the
remaining vertices, contradicting $\omega (f)\leq 2n+1$. Therefore $f(c)=1$.
Since a weight of at least $2n$ is required to totally Roman dominate the
remaining disjoint union of $n$ triangles, we conclude that $\gamma
_{tR}(D_{n})=2n+1$.

Now, let $f$ be any $\gamma _{tR}$-function on $D_{n}$. Then $\omega
(f)=2n+1 $ and $f(c)=1$. To dominate each triangle of $D_{n}-c$ with a
weight of $2 $, $\{f(u_{i}),f(v_{i})\}=\{0,2\}$ and $f(w_{i})=0$ for each $%
1\leq i\leq n $. Hence each $w_{i}$ is a dead vertex.~$\square $

{\begin{figure}
	\centering
	\begin{tikzpicture}		

	\begin{scope}[shift={(-3,0.75)}] 
		\node [std] (c) at (0, 0)[label=above: $c$]{};

		\foreach \i in {1,2,3}
			{\node [std] (u\i) at (120*\i:13mm)[label=above: $u_{\i}$] {};
			\node [std] (v\i) at (60+120*\i:13mm)[label=above: $v_{\i}$] {};
			\node [std] (w\i) at (30+120*\i:24mm)[label=above: $w_{\i}$] {};

			\draw (c)--(u\i);
			\draw (c)--(v\i);
			\draw (w\i)--(u\i);
			\draw (w\i)--(v\i);
			\draw (u\i)--(v\i);}
	\end{scope}

	\begin{scope}[shift={(3,0)}] 
		\node [std] (c') at (0, 0)[label=above: $c$]{};

		\foreach \i in {1,2,3,4}
			{\node [std] (u\i') at (22.5+90*\i:15mm)[label=right: $u_{\i}$] {};
			\node [std] (v\i') at (67.5+90*\i:15mm)[label=right: $v_{\i}$] {};
			\node [std] (w\i') at (45+90*\i:27mm)[label=right: $w_{\i}$] {};

			\draw (c')--(u\i');
			\draw (c')--(v\i');
			\draw (w\i')--(u\i');
			\draw (w\i')--(v\i');
			\draw (u\i')--(v\i');}
	\end{scope}

	\end{tikzpicture}	
	\caption{The graphs $D_3$ and $D_4$}
\end{figure}
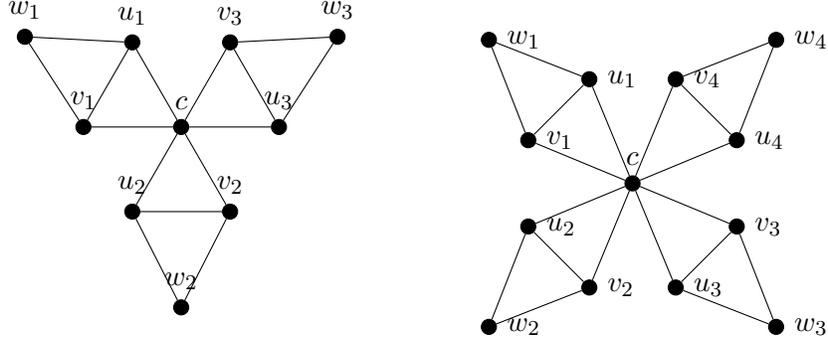}

\vspace{1.2mm}

The following result shows that, for $n\geq 3$, every edge in $E(\overline{%
D_{n}})$ incident with $w_{i}$ is critical.

\begin{prop}
If $n\geq 3$, $i\in \{1,...,n\}$ and $w_{i}v\in E(\overline{D_{n}})$, then $%
\gamma _{tR}(D_{n}+w_{i}v)<\gamma _{tR}(D_{n})$.
\end{prop}

\noindent \emph{Proof.} Without loss of generality, consider an edge $%
w_{1}v\in E(\overline{D_{n}})$. Then (without loss of generality) $v\in
\{w_{2},u_{2},c\}$. If $v=w_{2}$, define $f:V(D_{n}+w_{1}v)\rightarrow
\{0,1,2\}$ by $f(w_{1})=f(w_{2})=1$, $f(c)=f(u_{3})=\cdots =f(u_{n})=2$, and 
$f(y)=0$ for all other $y\in V(D_{n})$. Otherwise, if $v\in \{u_{2},c\}$,
define $f:V(D_{n}+w_{1}v)\rightarrow \{0,1,2\}$ by $f(c)=f(u_{2})=f(u_{3})=%
\cdots =f(u_{n})=2$ and $f(y)=0$ for all other $y\in V(D_{n})$. In either
case, $f$ is a TRD-function on $D_{n}+w_{1}v$ and $\omega (f)=2n$.
Therefore, by \thref{tR(D_n)}, every edge $w_{i}v\in E(\overline{D_{n}})$ is
critical.~$\square $

\vspace{1.2mm}

However, for $n\geq 3$, $D_{n}$ is not $\gamma _{tR}$-edge-critical since
(for example) $\gamma _{tR}(D_{n}+u_{1}u_{2})=2n+1$. Furthermore, $D_{2}$ is
not $\gamma _{tR}$-edge-critical since (for example) $\gamma
_{tR}(D_{2}+w_{1}w_{2})=5$. However, adding edges to $D_{n}$ until a $(2n+1)$%
-$\gamma _{tR}$-edge-critical graph $D_{n}^{\prime }$ is obtained results in 
$D_{n}^{\prime }$ having no dead vertices. Hence we pose the following
question.

\begin{que}
Do there exist $\gamma _{tR}$-edge-critical graphs containing dead vertices?
\end{que}

\vspace{1.2mm}

We characterized $\gamma _{tR}$-edge-critical spiders in 
\thref{edge-crit spider}. Finding other classes of $\gamma _{tR}$-edge-critical trees and,
indeed, characterizing $\gamma _{tR}$-edge-critical trees, remain open
problems.

\pagebreak

\label{Refs}

\end{document}